\documentclass[11pt]{article}
\usepackage{graphicx}
\usepackage{epsfig}
\usepackage{latexsym}
\newsavebox{\toy}
\savebox{\toy}{\framebox[0.65em]{\rule{0cm}{1ex}}}
\newcommand{\QED}{\usebox{\toy}\end{demo}}
\newenvironment{property}%
{\begin{list}{}{\setlength{\rightmargin}{0pt}%
\setlength{\itemsep}{0pt}}}{\end{list}}
\newlength{\templength}
\newcommand{\bp}{\setlength{\templength}{\labelwidth}%
\setlength{\labelwidth}{2em}\begin{property}}
\newcommand{\ep}{\end{property}\setlength{\labelwidth}{\templength}}
\newtheorem{theorem}{Theorem}[subsection]
\newtheorem{lemma}[theorem]{Lemma}
\newtheorem{proposition}[theorem]{Proposition}
\newtheorem{corollary}[theorem]{Corollary}
\newtheorem{assumption}{Assumption}
\newtheorem{definition}{Definition}[subsection]
\newtheorem{remark}{Remark}[subsection]
\newtheorem{exercise}{Exercise}[subsection]

\newcommand{\Thm}[1]{Theorem \ref{Thm.#1}}
\newcommand{\Lem}[1]{Lemma \ref{Lem.#1}}

\newcommand{\Prop}[1]{Proposition \ref{Prop.#1}}

\newcommand{\Theorem}[1]{\begin{theorem}\label{Thm.#1}}
\newcommand{\Lemma}[1]{\begin{lemma}\label{Lem.#1}}
\newcommand{\Proposition}[1]{\begin{proposition}\label{Prop.#1}}
\newcommand{\Corollary}[1]{\begin{corollary}\label{Cor.#1}}
\newcommand{\Assumption}[1]{\begin{assumption}\label{Ass.#1}\rm}
\newcommand{\Definition}[1]{\begin{definition}\label{Def.#1}\rm}
\newcommand{\Remark}[1]{\begin{remark}\label{Rem.#1}\rm }
\newcommand{\Exercise}[1]{\begin{exercise}\label{Exe.#1}\rm }

\newcommand{\bd}{\begin{displaymath}}
\newcommand{\ed}{\end{displaymath}}
\newcommand{\bdn}{\begin{equation}}
\newcommand{\bdnl}{\begin{equation}\label}
\newcommand{\edn}{\end{equation}}
\newcommand{\barray}{\begin{array}}
\newcommand{\earray}{\end{array}}
\newcommand{\bds}{\begin{description}}
\newcommand{\eds}{\end{description}}
\newcommand{\bitemize}{\begin{itemize}}
\newcommand{\eitemize}{\end{itemize}}
\newcommand{\benumerate}{\begin{enumerate}}
\newcommand{\eenumerate}{\end{enumerate}}
\newcommand{\btabbing}{\begin{tabbing}}
\newcommand{\etabbing}{\end{tabbing}}
\newcommand{\bcenter}{\begin{center}}
\newcommand{\ecenter}{\end{center}}
\newcommand{\bflushright}{\begin{flushright}}
\newcommand{\bflushleft}{\begin{flushleft}}
\newcommand{\eflushright}{\end{flushright}}
\newcommand{\eflushleft}{\end{flushleft}}
\newcommand{\bdnn }{\begin{eqnarray*}}
\newcommand{\ednn }{\end{eqnarray*}}
\newcommand{\bdmn}{\begin{eqnarray}}
\newcommand{\edmn}{\end{eqnarray}}

\newcommand{\nn}{\nonumber}
\newcommand{\SSC}[1]{\section{#1}\setcounter{equation}{0}}

\newcounter{biblio}
\newenvironment{references}%
{\begin{list}{[\arabic{biblio}]}{\usecounter{biblio}%
\setlength{\leftmargin}{2.5em}\setlength{\rightmargin}{0pt}%
\setlength{\labelwidth}{2em}\setlength{\itemsep}{0pt}}}{\end{list}}
\newcommand{\References}%
{\vspace{2.8ex plus .3ex minus .3ex}%
\begin{center}{\bf References}\end{center}\begin{references}}

\usepackage{amssymb}
\usepackage{mathrsfs}
\newcommand{\bL}{{\mathbb{L}}}
\newcommand{\N}{{\mathbb{N}}}
\newcommand{\Z}{{\mathbb{Z}}}
\newcommand{\zd}{\Z^d}

\newcommand{\R}{{\mathbb{R}}}
\newcommand{\rd}{\R^d}

\newcommand{\ra }{\rightarrow }
\newcommand{\lra }{\longrightarrow }

\newcommand{\Ra}{\Rightarrow }

\newcommand{\ov}{\overline}
\newcommand{\tl}{\widetilde}

\newcommand{\vvs}{\vspace{2ex}}
\newcommand{\vs}{\vspace{1ex}}

\newcommand{\lan}{\langle \:}
\newcommand{\ran}{\: \rangle}
\newcommand{\lef}{\left}
\newcommand{\rig}{\right}
\newcommand{\ri}{\right}
\newcommand{\st}{\stackrel}

\newcommand{\8}{\infty}

\newcommand{\dps}{\displaystyle}
\newcommand{\sub}{\subset}

\newcommand{\bsh}{\backslash}
\newcommand{\pri}{\prime}

\newcommand{\inflim}{\mathop{\underline{\lim}}}
\newcommand{\suplim}{\mathop{\overline{\lim}}}

\renewcommand{\a}{\alpha}
\renewcommand{\b}{\beta}
\newcommand{\gm}{\gamma}

\newcommand{\del}{\delta}
\newcommand{\D}{\Delta}

\newcommand{\e}{\varepsilon}

\newcommand{\z}{\zeta}
\newcommand{\h}{\eta}

\newcommand{\lm}{\lambda}

\newcommand{\rh}{\rho}

\newcommand{\s}{\sigma}

\newcommand{\W}{\Omega}

\newcommand{\cF }{{\cal F}}

\newcommand{\cO }{{\cal O}}

\newcommand{\cR }{{\cal R}}

\newcommand{\ovn}{\ov{N}}
\newcommand{\ovN}{\ov{N}}

\makeatletter
\def\section{\@startsection{section}{1}{\z@}{-3.5ex plus -1ex minus 
 -.2ex}{2.3ex plus .2ex}{\bf}}
\def\subsection{\@startsection{subsection}{2}{\z@}{-3.25ex plus -1ex minus 
 -.2ex}{1.5ex plus .2ex}{\bf}}
\makeatother

\begin{document}
\bcenter

\large{\bf Localization for   
Linear Stochastic Evolutions}\footnote{\today }\\

\vvs \normalsize

\noindent Nobuo YOSHIDA\footnote{ 
Division of Mathematics,
Graduate School of Science,
Kyoto University,
Kyoto 606-8502, Japan.
email: {\tt nobuo@math.kyoto-u.ac.jp}
URL: {\tt http://www.math.kyoto-u.ac.jp/}$\widetilde{}$ {\tt nobuo/}.
Supported in part by JSPS Grant-in-Aid for Scientific
Research, Kiban (C) 17540112}

\ecenter
%%%%%%%%%%%%%%%%%%%%
\begin{abstract}
We consider a discrete-time stochastic growth model on the 
$d$-dimensional lattice with non-negative real numbers 
as possible values per site. 
The growth model describes various interesting examples such as 
oriented site/bond percolation, directed polymers in random environment, 
time discretizations of the binary contact path process. 
We show the equivalence between the
slow population growth and a localization property in
terms of ``replica overlap". The main novelty of this paper is 
that we obtain this equivalence 
even for models with positive probability of extinction 
at finite time. 
In the course of the proof, we characterize, in a general setting, 
 the event on which an exponential martingale vanishes in the limit. 
\end{abstract}

%%%%%%%%%%%%%%%%%%%%
\small
\noindent AMS 2000 subject classification: Primary 60K35; 
secondary 60J37, 60K37, 82B26.\\
Key words and phrases: 
linear stochastic evolutions, localization, slow growth phase.

\tableofcontents

\normalsize
%%%%%%
\SSC{Introduction}
%%%%%%
We write $\N=\{0,1,2,...\}$, 
$\N^*=\{1,2,...\}$ and 
$\Z=\{ \pm x \; ; \; x \in \N \}$. For 
$x=(x_1,..,x_d) \in \rd$, $|x|$ stands for the $\ell^1$-norm: 
$|x|=\sum_{i=1}^d|x_i|$. For $\xi=(\xi_x)_{x \in \zd} \in \R^{\zd}$, 
$|\xi |=\sum_{x \in \zd}|\xi_x|$. 
Let 
$(\W, \cF, P)$ be 
a probability space.
We write $P[X]=\int X \; dP$ and 
$P[X:A]=\int_A X \; dP$ for a random variable $X$ and an 
event $A$. For events $A,B \sub \W$, $A \sub B$ a.s. means that 
$P(A \bsh B)=0$. Similarly, $A = B$ a.s. means that 
$P(A \bsh B)=P(B \bsh A)=0$.
%%%%%%%%%%%%%%%%
\subsection{The oriented site percolation (OSP)}
%%%%%%%%%%%%%
We start by discussing the {\it oriented site 
percolation} as a motivating example.
Let 
$\h_{t,y}$, $(t,y) \in \N^*\times \zd$ be $\{0,1\}$-valued i.i.d. 
random variables  
with $P(\h_{t,y}=1)=p \in (0,1)$. The site $(t,y)$ with 
$\h_{t,y}=1$ and $\h_{t,y}=0$ are referred to respectively as 
{\it open} and {\it closed}.
An {\it open oriented path} from $(0,0)$ to 
$(t,y) \in \N^*\times \zd$ is a sequence 
$\{(s,x_s) \}_{s=0}^t$ in $\N \times \zd$ such that 
$x_0=0$, $x_t=y$, 
$|x_s-x_{s-1}|=1$, 
$\h_{s,x_s}=1$ for all $s=1,..,t$.
For oriented percolation, 
it is traditional to discuss the presence/absence of the 
open oriented paths to certain time-space location. 
On the other hand, the model exhibits another type 
of phase transition, if we look at 
not only the presence/absence of the open oriented paths, 
but also their number.
Let $N_{t,y}$  be the number of 
open oriented paths from $(0,0)$ to 
$(t,y)$ and let $|N_t|=\sum_{y \in \zd}N_{t,y}$ be the 
total number of  open oriented paths from $(0,0)$ to the ``level'' $t$. 
If we regard each open oriented path 
$\{(s,x_s) \}_{s=0}^t$ as a trajectory of a particle, 
then $N_{t,y}$ is the number of the particles which 
occupy the site $y$ at time $t$. 

We now  note that 
$|\ov{N}_t|\st{\rm def.}{=}(2dp)^{-t}|N_t|$ is a martingale, 
since each open oriented path from $(0,0)$ to $(t,y)$ branches 
and survives to the next level via $2d$ neighbors of $y$, each of 
which is open with probability $p$. Thus, 
by the martingale convergence theorem, the following limit 
exists a.s.:
$$
|\ovn_\8|\st{\rm def}{=}\lim_{t \ra \8}|\ovn_t|.
$$
Moreover, 
\bds
\item[i)] 
If $d \ge 3$ and $p$ is large enough,  
then, $P(|\ovn_\8|>0)>0$, 
which means that, at least with positive probability, the 
total number of paths 
$|N_t|$ is of the same order as its expectation $(2pd)^t$ 
as $t \ra \8$.
\item[ii)] 
If $d =1,2$, then for all $p \in (0,1)$, $P(|\ovn_\8|=0)=1$, 
which means that the 
total number of paths $|N_t|$ is of smaller order 
than its expectation $(2pd)^t$ a.s. as $t \ra \8$.
Moreover, there is a non-random constant $c>0$ such that 
$|\ovn_t|=\cO (\exp (-ct))$ a.s. as $t \ra \8$.
\eds
This phase transition was predicted by T. Shiga in late 1990's and 
the proof was given recently in \cite{Ber09,Yo08b}.

We denote the
density of the population by: 
\bdnl{rh}
\rh_t(x)=\frac{N_{t,x}}{|N_t|}{\bf 1}_{\{ |N_t|>0\}}, 
\; \; t \in \N, x \in \zd.
\edn 
Here and in what follows, we adopt the following convention.
For a random variable $X$ defined on an event $A$, we define the 
random variable $X{\bf 1}_A$ by  $X{\bf 1}_A=X$ on $A$ and 
$X{\bf 1}_A=0$ outside $A$.
Interesting objects related to the
density would be 
\bdnl{rh^*} \rh^*_t=\max_{x \in \zd}\rh_t (x), \;
\; \mbox{and}\; \; \cR_t=|\rh_t^2|=\sum_{x \in \zd}\rh_t(x)^2. 
\edn
$\rh^*_t$ is the density at the most populated site, while $\cR_t$
is the probability that two particles picked up 
randomly from the total population at time $t$ are
at the same site. We call $\cR_t$ the {\it replica overlap}, in
analogy with the spin glass theory. Clearly, $(\rh^*_t)^{2} \le \cR_t
\le \rh^*_t$. These quantities convey information on
localization/delocalization of the particles. Roughly speaking,
large values of $\rh^*_t$ or $\cR_t$ indicate that most of
the particles are concentrated on small numbers of ``favorite
sites" ({\it localization}), whereas small values of them imply 
that the particles are spread out over large number of sites ({\it
delocalization}).

As applications of results in this paper, we get the following 
result. 
It says that, in the presence of an  
infinite open path, the slow growth $|\ovN_\8|=0$ is equivalent to
a localization property $\suplim_{t \ra \8}\cR_t  \ge c>0$.
Here, and in what follows, a {\it constant} always means a 
{\it non-random constant}. 
%%%%%%%%%%%%%%%%%%
\Theorem{OPloc}
%%%%%%%%%
\bds
\item[a)]
If $P(|\ovn_\8|>0)>0$, then, 
${\dps 
\sum_{t \ge 1}\cR_t<\8\; \; a.s.}$ 
\item[b)] If $P(|\ovn_\8|=0)=1$, 
then, there exists a constant $c >0$ such that:
\bdnl{OPlimR>c}
\{  \mbox{$|N_t|>0$ for all $t \in \N$}\} 
=
\lef\{ \suplim_{t \ra \8}\cR_t  \ge c \rig\}
\; \; \; \mbox{a.s.}
\edn
\eds
%%%%%%%%%%%%%%%
\end{theorem}
%%%%%%%%%
Note that $P(|\ovN_\8|=0)=1$ for all $p \in (0,1)$ if $d \le 2$. 
Thus, (\ref{OPlimR>c}) in particular means that, if $d \le 2$, 
the path localization $\suplim_{t \ra \8}\cR_t  \ge c$ occurs 
a.s. on the event of percolation. 
\Thm{OPloc} is shown at the end of section \ref{sec:results} 
as a consequence of more general results for linear 
stochastic evolutions.
%%%%%%%
\subsection{The linear stochastic evolution} \label{sec:lse}
%%%%%%%%%
We now introduce the framework of this article.
Let $A_t =(A_{t,x,y})_{x,y \in \zd}$, $t \in \N^*$ 
be a sequence of random matrices on a probability space 
$(\W, \cF, P)$ such that:
\bdn
\mbox{ $A_1,A_2,...$ are i.i.d.} \label{Aiid}
\edn
Here are the set of assumptions we assume for $A_1$:
\bdmn 
& & \mbox{ $A_1$ is not a constant matrix.} \label{AnotC} \\
& & \mbox{ $A_{1,x,y} \ge 0$ for all $x,y \in \zd$.} \label{A>0} \\
& & \mbox{The columns 
$\{ A_{1,\cdot,y} \}_{y \in \zd}$ are independent.} 
\label{colind}\\
& & P[A_{1,x,y}^3]<\8
\; \; \; \mbox{for all $x,y \in \zd$,}
\label{A^3} \\
& & A_{1,x,y}=0 \; \; \mbox{a.s. if $|x-y| >r_A$ for some 
non-random $r_A \in \N$.} \label{r_A} \\
& & 
\mbox{$(A_{1,x+z,y+z})_{x,y \in \zd}
\st{\rm law}{=}A_1$ for all $z \in \zd$}.\label{A=A} \\
& & \barray{l}
\mbox{The set $ \{x \in \zd\; ; \; \sum_{y \in \zd}a_{x+y}a_y \neq 0\}$ 
contains a linear basis of $\rd$,}\\
\mbox{where $a_y=P[A_{1,0,y}]$.}
\earray \label{irred}
\edmn 
Depending on the results we prove in the sequel, some of 
these conditions can be relaxed. However, 
we choose not to bother ourselves with the pursuit of 
the minimum assumptions for each result.

We define a Markov chain 
$(N_t)_{t \in \N}$ with values in $[0,\8)^{\zd}$ by:
\bdnl{lse}
\sum_{x \in \zd}N_{t-1,x}A_{t,x,y}=N_{t,y}, \; \; t \in \N^*.
\edn 
In this article, we suppose that the initial state $N_0$ is given by 
``a single particle at the origin":
\bdnl{N_0}
N_0=\lef( \del_{0,x}\ri)_{x \in \zd}
\edn
Here and in what follows, $\del_{x,y}={\bf 1}_{\{x=y\}}$ for 
$x,y \in \zd$.  
If we regard $N_t \in [0,\8)^{\zd}$ as a 
row vector, (\ref{lse}) can be interpreted as: 
$$
N_t=N_0A_1A_2\cdots A_t, \; \; \; t=1,2,...
$$
The Markov chain defined above 
can be thought of as the time discretization of the linear particle 
system considered in the last Chapter in T. Liggett's book 
\cite[Chapter IX]{Lig85}. Thanks to the time discretization, 
the definition is considerably simpler here. 
Though we {\it do not} assume in general that $(N_t)_{t \in \N}$ takes 
values in $\N^{\zd}$, we refer $N_{t,y}$ as the 
``number of particles" at time-space $(t,y)$, and 
$|N_t|$ as the ``total number of particles" at time $t$.

We now see that various interesting examples 
are included in this framework. We recall the notation 
$a_y$ from (\ref{irred}).

\vvs
\noindent $\bullet$
%%%%%%%%%%%%%%%%%%%%%%%%
{\bf Generalized oriented site percolation (GOSP):}
%%%%%%%%%%%%%%%%%%%%%%%%%%%%%
We generalize OSP as follows. 
Let $\h_{t,y}$, $(t,y) \in \N^*\times \zd$ be $\{0,1\}$-valued i.i.d. 
random variables with $P(\h_{t,y}=1)=p \in [0,1]$ and let 
$\z_{t,y}$, $(t,y) \in \N^*\times \zd$ be 
another $\{0,1\}$-valued i.i.d. random variables  with 
$P(\z_{t,y} =1)=q \in [0,1]$, which are independent of $\h_{t,y}$'s. 
To exclude trivialities, we assume that either $p$ or $q$ is in 
$(0,1)$. 
We refer to the process $(N_t)_{t \in \N}$ defined by (\ref{lse}) with:
$$
A_{t,x,y}={\bf 1}_{|x-y|=1}\h_{t,y}+\del_{x,y}\z_{t,y}
$$
as the {\it generalized oriented site percolation} (GOSP).
Thus, the OSP is the special case ($q=0$) of GOSP.
%%%%%%%%%
%The GOP can be interpreted as 
%a time discretization of the binary contact path process \cite{Gri83}, 
%\cite[Chapter IX]{Lig85}, \cite{NaYo08}.
%%%%%%%%
The covariances of $(A_{t,x,y})_{x,y \in \zd}$ can be 
seen from:
\bdnl{OP}
a_y = p{\bf 1}_{\{|y|=1\}}+q\del_{y,0}, \; \; \; 
P[A_{t,x,y}A_{t,\tl{x},y}]=
\lef\{ \barray{ll}
q & \mbox{if $x=\tl{x}=y$,}\\
p & \mbox{if $|x-y|=|\tl{x}-y|=1$,}\\
a_{y-x}a_{y-\tl{x}} & \mbox{if otherwise.}
\earray  \ri.
\edn
In particular, we have $|a|=2dp+q$ (Recall that $|a|=\sum_ya_y$). 
%%%%%%%%
%\bdnl{OP}
%a_y = p{\bf 1}_{|y|=1}, \; \; \; 
%P[A_{t,x,y}A_{t,\tl{x},\tl{y}}]=
%\lef\{ \barray{ll}
%a_{y-x}a_{\tl{y}-\tl{x}} & \mbox{if $y \neq \tl{y}$,}\\
%a_{y-x}a_{y-\tl{x}}p^{-1} & \mbox{if $y=\tl{y}$.}
%\earray  \ri.
%\edn
%%%%%%%%%
%%%%%%%%
\vs

%%%%%%%%%%%%%%%%%%%%%%%%
\noindent $\bullet$ {\bf Generalized oriented bond percolation (GOBP):}
%%%%%%%%%%%%%%%%%%%%%%%%%%%%%
Let $\h_{t,x,y}$, $(t,x,y) \in \N^*\times \zd \times \zd$ 
be $\{0,1\}$-valued i.i.d.random variables 
with $P(\h_{t,x,y}=1)=p \in [0,1]$ and let 
$\z_{t,y}$, $(t,y) \in \N^*\times \zd$ be 
another $\{0,1\}$-valued i.i.d. random variables with 
$P(\z_{t,y} =1)=q \in [0,1]$, which are independent of $\h_{t,y}$'s.
We refer to the process $(N_t)_{t \in \N}$ defined by (\ref{lse}) with:
$$
A_{t,x,y}={\bf 1}_{\{|x-y|=1\}}\h_{t,x,y}+\del_{x,y}\z_{t,y}
$$
as the {\it generalized oriented bond percolation} (GOBP). 
We call the special case $q=0$ {\it oriented bond percolation} (OBP).
To interpret the definition, let us call the pair of time-space points 
$\lan (t-1,x), (t,y) \ran$ 
a {\it bond} if $|x-y| \le 1$, $(t,x,y) \in \N^*\times \zd \times \zd$.
A bond $\lan (t-1,x), (t,y) \ran$ with $|x-y| = 1$ is said 
to be {\it open} if $\h_{t,x,y}=1$, and a bond 
$\lan (t-1,y), (t,y) \ran$ is said 
to be {\it open} if $\z_{t,y}=1$. 
For GOBP, an {\it open oriented path} from $(0,0)$ to 
$(t,y) \in \N^*\times \zd$ is a sequence 
$\{(s,x_s) \}_{s=0}^t$ in $\N \times \zd$ such that 
$x_0=0$, $x_t=y$ and bonds 
$\lan (s-1,x_{s-1}), (s,x_s) \ran$ are open 
for all $s=1,..,t$. If $N_0=(\del_{0,y})_{y \in \zd}$, then, 
the number of open oriented paths from 
$(0,0)$ to $(t,y) \in \N^*\times \zd$ is given by $N_{t,y}$. 

The covariances of $(A_{t,x,y})_{x,y \in \zd}$ can be 
seen from:
\bdnl{OBP}
a_y = p{\bf 1}_{\{|y|=1\}}+q\del_{y,0}, \; \; \; 
P[A_{t,x,y}A_{t,\tl{x},y}]=
\lef\{ \barray{ll}
a_{y-x} & \mbox{if $x=\tl{x}$,}\\
%%%%%%%%%%%%%%%%%%%%%%%%%%%%%%%%%%%%%%%%%%%%%%%%%%%%%%
%pq & \mbox{if $y=\tl{y}$, $x=y$, $|\tl{x}-y|=1$,} \\
%pq & \mbox{if $y=\tl{y}$, $|x-y|=1$, $\tl{x}=y$,} \\
%%%%%%%%%%%%%%%%%%%%%%%%%%%%%%%%%%%%%%%%%%%%%%%%%%%%%%%
a_{y-x}a_{y-\tl{x}} & \mbox{if otherwise.}
\earray  \ri.
\edn
In particular, we have $|a|=2dp+q$. 
%%%%%%%%%%%%%%%%%%%%%%%%%%%%%%%%%%%%%%

\vs
\noindent $\bullet$
%%%%%%%%%%%%%%%%%%%%%%%%
{\bf Directed polymers in random environment (DPRE):}
%%%%%%%%%%%%%%%%%%%%%%%%%%%%%
Let $\{\h_{t,y} \; ; \; (t,y) \in \N^*\times \zd \}$ be i.i.d. 
with $\exp (\lm (\b))\st{\rm def.}{=}
P[\exp (\b \h_{t,y})]<\8$ for any $\b \in (0,\8)$.
The following expectation is called the 
partition function of the {\it directed polymers in random environment}:
$$
N_{t,y}=P^0_S\lef[ \exp \lef( \b\sum_{u=1}^t\h_{u,S_u}\ri):S_t=y\ri], 
\; \; \; (t,y) \in \N^* \times \zd,
$$
where $((S_t)_{t \in \N},P_S^x)$ is the simple random walk on $\zd$.
We refer the reader to a review paper 
\cite{CSY04} and the references therein for more 
information. 
Starting from $N_0=(\del_{0,x})_{x \in \zd}$, 
the above expectation can be obtained inductively by (\ref{lse}) 
with:
$$
A_{t,x,y}={{\bf 1}_{|x-y|=1} \over 2d}\exp (\b \h_{t,y}).
$$
The covariances of $(A_{t,x,y})_{x,y \in \zd}$ can be 
seen from:
\bdnl{DP}
a_y = {e^{\lm (\b)}{\bf 1}_{\{|y|=1\}} \over 2d}, \; \; \; 
P[A_{t,x,y}A_{t,\tl{x},y}]=
e^{\lm (2\b)-2\lm (\b)}
a_{y-x}a_{y-\tl{x}} 
\edn
In particular, we have $|a|=e^{\lm (\b)}$.
%%%%%%%%%

\vs
\noindent $\bullet$
%%%%%%%%%%%%%%%%%%%%%%%%
%%%%%%%%%%%%%%%%%%%%%%%%
{\bf The binary contact path process (BCPP):}
%%%%%%%%%%%%%%%%%%%%%%%%%%%%%
The binary contact path process is a continuous-time Markov 
process with values in  $\N^{\zd}$, originally introduced by 
D. Griffeath \cite{Gri83}. 
In this article, we consider a discrete-time variant as follows.
Let 
\bdnn
& & \{\h_{t,y}=0,1\; ; \; (t,y) \in \N^*\times \zd \}, 
\; \; \; \{\z_{t,y}=0,1\; ; \; (t,y) \in \N^*\times \zd \}, \\
& & \{e_{t,y}\; ; \; (t,y) \in \N^*\times \zd \}
\ednn
be families of 
i.i.d. random variables with $P(\h_{t,y}=1)=p \in (0,1]$, 
$P(\z_{t,y}=1)=q \in [0,1]$, 
and $P(e_{t,y}=e)={1 \over 2d}$ for each 
$e \in \zd$ with $|e|=1$. We suppose that these three 
families are independent of each other. Starting from 
an $N_0 \in \N^{\zd}$, we define a Markov chain 
$(N_t)_{t \in \N}$ with values in $\N^{\zd}$ by:
$$
N_{t+1,y}=\h_{t+1,y}N_{t,y-e_{t+1,y}}+\z_{t+1,y}N_{t,y}, \; \; \; 
t \in \N.
$$
We interpret the process as the spread of an infection, 
with $N_{t,y}$ infected individuals at time $t$ at the site $y$. 
The $\z_{t+1,y}N_{t,y}$ term above means that these 
individuals remain infected at time $t+1$ with probability $q$, 
and they recover with probability $1-q$. On the other hand, 
the $\h_{t+1,y}N_{t,y-e_{t+1,y}}$ term means that, 
with probability $p$, a neighboring site $y-e_{t+1,y}$ is 
picked at random (say, the wind blows from that direction), and 
$N_{t,y-e_{t+1,y}}$ individuals at site $y$ are infected anew 
at time $t+1$. 
This Markov chain is obtained by (\ref{lse}) with:
$$
A_{t,x,y}=\h_{t,y}{\bf 1}_{e_{t,y}=y-x}+\z_{t,y}\del_{x,y}.
$$
The covariances of $(A_{t,x,y})_{x,y \in \zd}$ can be 
seen from:
\bdnl{BCPP[AA]}
a_y  =  {p{\bf 1}_{\{|y|=1\}} \over 2d}+q\del_{0,y},
\; \; \; 
P[A_{t,x,y}A_{t,\tl{x},y}] = 
\lef\{ \barray{ll}
a_{y-x} & \mbox{if $x=\tl{x}$,}\\
\del_{x,y}qa_{y-\tl{x}}+
\del_{\tl{x},y}qa_{y-x}
& \mbox{if $x \neq \tl{x}$.}
\earray  \ri.
\edn
%%%%%%
%\bdmn
%a_y & = & {p{\bf 1}_{\{|y|=1\}} \over 2d}+q\del_{0,y},
%\label{BCPPa}   \\
%P[A_{t,x,y}A_{t,\tl{x},y}]& = &
%\lef\{ \barray{ll}
%a_{y-x} & \mbox{if $x=\tl{x}$,}\\
%q\del_{x,y}a_{y-\tl{x}}+
%q\del_{\tl{x},y}a_{y-x}
%& \mbox{if $x \neq \tl{x}$.}
%\earray  \ri.
%\label{BCPP[AA]} 
%\edmn
%%%%%%%%%%%
In particular, we have $|a|=p+q$. 
%%%%%%%%%

\vvs
\noindent {\bf Remark:} 
The branching random walk in random environment 
considered in \cite{HuYo08,Shi09a,Shi09b,Yo08a} can also be 
considered as a ``close relative" to the models considered here,
although it does not exactly fall into our framework.
%%%%%%%%%
\subsection{The regular and slow growth phases}
%%%%%%%%

We now recall the following facts and notion from 
\cite[Lemmas 1.3.1 and 1.3.2]{Yo08b}. Let $\cF_t$ be the 
$\s$-field generated by $A_1,..,A_t$.
%%%%%%%%%%%%%%
\Lemma{0,1}
%%%%%%%%%%
Define $\ov{N}_t=\lef( \ovn_{t,x}\ri)_{x \in \zd}$ by:
\bdnl{ovn_t}
\ovn_{t,x}=|a|^{-t}N_{t,x}.
\edn
\bds
\item[a)]
$(|\ov{N}_t|, \cF_t)_{t \in \N}$ is a martingale, 
and therefore, the following limit exists a.s. 
\bdnl{ovn_8}
|\ovn_\8| =\lim_{t \ra \8}|\ovn_t|.
\edn
\item[b)] Either 
\bdnl{0,1}
P[|\ovn_\8|]=1\; \; \mbox{or}\; \; 0.
\edn
Moreover, $P[|\ovn_\8|]=1$ if and only if the limit (\ref{ovn_8}) 
is convergent in $\bL^1 (P)$. 
\eds
\end{lemma}
%%%%%%%%%%
We will refer to the former case of (\ref{0,1})
as {\it regular growth phase} and the latter as 
{\it slow growth phase}. 

The regular growth 
means that, at least with positive probability, 
the growth of the ``total number" $|N_t|$ 
of particles is of the same order as its expectation $|a|^t|N_0|$.
On the other hand, the slow growth means that, almost surely, 
 the growth of $|N_t|$ is slower than its expectation. 

\vvs
We now recall from \cite{Ber09} and 
\cite[Theorems 3.1.1 and 3.2.1]{Yo08b} 
the following criterion for slow growth phase.
%%%%%%%%%%%%%%%%%%%%%%
\Proposition{r/s}
%%%%%%%%%
$P (|\ovn_\8|=0)=1$ if $d=1,2$, or if:
\bdnl{SGlog}
\sum_{y \in \zd}P\lef[ A_{1,0,y}\ln A_{1,0,y}\ri] >|a| \ln |a|.
\edn
%%%%%%%%
\end{proposition}
%%%%%%%%%
The condition (\ref{SGlog}) roughly says that the matrix 
$A_1$ is ``random enough''. For DPRE, (\ref{SGlog}) is equivalent to 
$\b \lm^\pri (\b)-\lm (\b)>\ln (2d)$. 
%%%%%%%%
\subsection{The results} \label{sec:results}
%%%%%%%%
We introduce the following additional condition, 
which says  that 
the entries of the matrix $A_1$ are 
positively correlated in the following weak sense:
there is a constant $\gm \in (1,\8)$ such that:
\bdnl{covA} 
\sum_{x,\tl{x},y \in \zd}
\lef( P[A_{1,x,y}A_{1, \tl{x},y}] -\gm a_{y-x}a_{y-\tl{x}} \ri)
\xi_x\xi_{\tl{x}}\ge  0
\edn
for all $\xi \in [0,\8)^{\zd}$ such that $|\xi|<\8$.
%%%%%%%%%%%%%

\vvs
\noindent {\bf Remark:}
Clearly, (\ref{covA}) is satisfied if
there is a constant $\gm \in (1,\8)$ such that:
\bdnl{covA*} 
P[A_{1,x,y}, A_{1, \tl{x},y}]  \ge \gm a_{y-x}a_{y-\tl{x}}
\; \; \; \mbox{for all $x,\tl{x},y \in \zd$.}
\edn
For OSP and DPRE, we see from 
(\ref{OP}) and (\ref{DP}) that (\ref{covA*}) holds with: 
$$
\gm=1/p \; \; \mbox{and}\; \; \exp(\lm (2 \b)-2\lm (\b))
$$  
respectively for OSP and DPRE. 
For GOSP, GOBP and BCPP, (\ref{covA*}) is no longer true. 
However, one can check (\ref{covA}) for them with: 
$$
\gm =1 +\lef\{ \barray{ll}
{2dp(1-p)+q(1-q) \over (2dp+q)^2} & \mbox{for GOSP and GOBP}, \\
{p(1-p)+q(1-q) \over (p+q)^2} & \mbox{for BCPP}
\earray \rig.
$$
\cite[Remarks after Theorem 3.2.1]{Yo08b}.

We define the density $\rh_t (x)$ and the replica overlap 
$\cR_t$ in the same way as (\ref{rh}) and (\ref{rh^*}).

We first show that, on the event of survival, 
the slow growth is equivalent to the localization:
%%%%%%%%
\Theorem{loc}
%%%%%%%
Suppose (\ref{covA}). 
\bds
\item[a)]
If $P(|\ovN_\8|>0)>0$, then
${\dps 
\sum_{t \ge 0}\cR_t <\8\; \; \mbox{a.s.}}$
\item[b)] If $P(|\ovN_\8|=0)=1$, then
\bdnl{S=L}
\{ \mbox{\rm survival}\} 
=
\lef\{ \sum_{t \ge 0}\cR_t =\8 \rig\}
\; \; \; \mbox{a.s.}
\edn
where $\{ \mbox{\rm survival}\}\st{\rm def}{=}
\{  \mbox{$|N_t|>0$ for all $t \in \N$}\}$. 
Moreover, 
there exists a constant $c>0$ such that almost surely,
\bdnl{M<}
|\ovN_t| \le  \exp \lef( -c
\sum_{1 \le s \le t-1}\cR_s\ri) 
\; \; \mbox{for all large enough $t$'s}
\edn
\eds
%%%%%%%
\end{theorem}
%%%%%%%%%
\noindent {\bf Remark:} As can be seen from the proof 
(cf. \Prop{locb}a) below), 
(\ref{S=L}) is true even without assuming (\ref{covA}) and with 
(\ref{A^3}) replaced by a weaker assumption:
\bdnl{A^2}
P[A_{1,x,y}^2]<\8
\; \; \; \mbox{for all $x,y \in \zd$.}
\edn
\Thm{loc} says that, conditionally on survival, 
the slow growth $|\ovN_\8|=0$ is 
equivalent to the localization $\sum_{t \ge 0}\cR_t =\8$. 
We emphasize that this is the first case in which a 
result of this type is obtained for 
models with positive probability of extinction 
at finite time (i.e.,$P(|N_t|=0)>0$ for finite $t$). 
Similar results have been known before only 
in the case where no extinction at finite time is allowed, i.e.,
  $|N_t|>0$ for all $t \ge 0$, e.g., 
\cite[Theorem 1.1]{CaHu02}, \cite[Theorem 1.1]{CSY03}, 
\cite[Theorem 2.3.2]{CY05}, \cite[Theorem 1.3.1]{HuYo08}. 
The argument in the previous literature is roughly to show that
\bdnl{-lnN_t}
-\ln |\ovn_t| \asymp \sum^{t-1}_{u=0}\cR_u\; \; \mbox{a.s. as $t \ra \8$}
\edn
by using Doob's decomposition of the supermartingale $\ln |\ovn_t|$ 
(``$\asymp$" above means the asymptotic upper and lower bounds with 
positive multiplicative constants). 
This argument does not seem to be directly 
transportable to the case where the total population may get 
extinct at finite time, since $\ln |\ovn_t|$  is not even defined. 
To cope with this problem, we first characterize, in a general setting, 
the event on which an exponential martingale vanishes in the limit 
(\Prop{gloc} below). We then apply this characterization to the martingale 
$|\ovn_t|$. See also \cite{NY09b} for the application of this idea to 
the continuous-time setting.
%%%%%%%%%%%%%%%%%%%%%%%%%%%

\vvs
Next, we present a result which says that, under a mild assumption, 
we can replace 
$$
\sum_{t \ge 0}\cR_t =\8
$$
 in (\ref{S=L}) by 
a stronger localization property: 
$$
\suplim_{t \ra \8}\cR_t  \ge c, 
$$
where $c>0$ is a constant. 
To state the theorem,  we introduce 
some notation related to the random walk associated to our model.
Let $((S_t)_{t \in \N}, P_S^x)$ be the random walk on $\zd$ such that: 
\bdnl{S_t}
\mbox{$P_S^x(S_0=x)=1$ and $P_S^x(S_1=y)=a_{y-x}/|a|$}
\edn
and let $(\tl{S}_t)_{t \in \N}$ be its independent copy.
We then define: 
\bdnl{pi_d}
\pi_d =P^0_S \otimes P^0_{\tl{S}} ( \mbox{$S_t =\tl{S}_t$ for some $t \ge 1$}).
\edn
Then, by (\ref{irred}),
\bdnl{pi_d=1}
  \mbox{$\pi_d=1$ for $d =1,2$ and $\pi_d<1$ for $d \ge 3$}
\edn
%%%%%%%%
\Theorem{sloc}
%%%%%%%%
Suppose (\ref{covA}) and either of
\bds
\item[a)] $d=1,2$,
\item[b)] $P(|\ovN_\8|=0)=1$ and 
\bdnl{P(S=S)}
\gm >{ 1 \over \pi_d}, 
\edn
where $\gm$ and $\pi_d$ are from (\ref{covA}) and (\ref{pi_d}).
\eds
Then, there exists a constant $c>0$ such that:
\bdnl{limR>c}
\lef\{ \mbox{\rm survival} \rig\}
=\lef\{ \suplim_{t \ra \8}\cR_t  \ge c \rig\}
\; \; \; \mbox{a.s.}
\edn
%%%%%%%
\end{theorem}
%%%%%%%%%
This result generalizes \cite[Theorem 1.2]{CaHu02} 
and \cite[Proposition 1.4 b)]{CSY03}, which are 
obtained in the context of DPRE.
Similar results are also known for 
branching random walk in random environment 
\cite[Theorem 1.3.2]{HuYo08}. 
To prove \Thm{sloc}, we will use the argument which was 
initially applied to DPRE by P. Carmona and Y. Hu 
in \cite{CaHu02} (See also \cite{HuYo08}).
What is new in the present paper is to carry the arguments in the 
above mentioned papers over to the case where the extinction at 
finite time is possible. This will be done in section \ref{p:sloc1}.

\vvs
\noindent {\bf Remarks}
%%%%%%%%%%%%%
%%%%%%%%%%%
\noindent {\bf 1)}
%%%%%%%%
We prove (\ref{limR>c}) by way of the following stronger estimate:
$$
\inflim_{t \nearrow \8} {\sum_{s=0}^t\cR_s^{3/2}
\over \sum_{s=0}^t\cR_s} \ge c_1, \; \; a.s. 
$$ 
for some constant  $c_1>0$. This in particular implies the following
quantitative lower bound on the number of times at which the
replica overlap is larger than a certain positive number:
$$
\inflim_{t \nearrow \8} {\sum_{s=0}^t1_{\{ \cR_s \ge c_2 \}}
\over \sum_{s=0}^t\cR_s}
 \ge c_3, \; \; a.s.
$$
where $c_2$ and $c_3$ are positive constants 
(The inequality $r^{3/2} 
\le {\bf 1}\{ r \ge c\}+\sqrt{c}r$ for $r,c \in [0,1]$ can be used here). \\
%%%%%%%%%%%%%%%%%%%%%%%%%%%%%%%%%%%%%%%%%%%%%%%%%%%%%%%%%%%%%%
%%%%%%%%%%%
\noindent {\bf 2)}
%%%%%%%%
(\ref{limR>c}) is in contrast with the following 
delocalization result by M. Nakashima \cite{Nak09}:
if $d \ge 3$ and $\sup_{t \ge 0}P[|\ovn_t|^2]<\8$, then,
$$
\cR_t=\cO (t^{-d/2})\; \; 
\mbox{in $P (\; \cdot \; | |\ovn_\8|>0)$-probability .}
$$
See also \cite{NaYo08} for the continuous-time case and \cite{Shi09a,Yo08a} 
for the case of branching random walk in random environment.

%%%%%%%%%%
\vvs
Finally, we state the following variant of \Thm{sloc}, which says that even 
for $d \ge 3$, 
(\ref{P(S=S)}) can be dropped at the cost of some alternative assumptions.
Following M. Birkner \cite[page 81, (5.1)]{Bir03}, 
we introduce the following condition:
\bdnl{Cor4}
\sup_{t \in \N, x \in \zd}
{P_S^0 (S_t=x) \over P_S^0 \otimes P_{\tl{S}}^0(S_t=\tl{S}_t)}<\8,
\edn
which  is obviously 
true for the symmetric simple random walk on $\zd$.
%%%%%%%%
\Theorem{sloc2}
%%%%%%%%
Suppose $d \ge 3$, (\ref{covA}), (\ref{Cor4}) 
and that there exist mean-one i.i.d. random variables 
$\ov{\h}_{t,y}$, $(t,y) \in \N \times \zd$ such that:
\bdnl{ovh1}
A_{t,x,y}=\ov{\h}_{t,y}a_{y-x}.
\edn
Then, the slow growth ($P(|N_\8|=0)=1$) implies that 
there exists a constant $c >0$ such that (\ref{limR>c}) holds.
%%%%%%%
\end{theorem}
%%%%%%%%%
Note that OSP and DPRE  for $d \ge 3$ satisfy 
all the assumptions for \Thm{sloc2}. The proof of 
\Thm{sloc2} is based on \Thm{sloc} and a criterion 
for the regular growth phase, which is essentially due to 
M. Birkner \cite{Bir04}. These will be explained in 
section \ref{p:sloc2}.

\vvs
%%%%%%%%%
\noindent {\bf Proof of \Thm{OPloc}}:
%%%%%%%%%
The theorem follows from  \Thm{loc} and \Thm{sloc2}.
\hfill $\Box$.
%%%%%%%%%%%%%%
\SSC{Proofs of \Thm{loc}} \label{p:loc}
%%%%%%%%%
We will prove part b) first, and then part a).
%%%%%%%
\subsection{An abstraction of \Thm{loc}b)}
%%%%%%%
We will prove \Thm{loc}b) in the following generalized form, 
where the slow glowth ($P(|\ov{N}_\8|=0)=1$) is not assumed in advance:
%%%%%%%%%%
\Proposition{locb}
%%%%%%%
\bds
\item[a)] 
Even without assuming (\ref{covA}) and with 
(\ref{A^3}) replaced by (\ref{A^2}), it holds that
\bdnl{SsupL}
\{ |\ov{N}_\8| >0\} 
\supset 
\lef\{ \mbox{\rm survival}, \; \sum_{t \ge 0}\cR_t <\8 \rig\}
\; \; \; \mbox{a.s.}
\edn
\item[b)] Suppose (\ref{A^3}) and (\ref{covA}).
Then, there exists a constant $c>0$ such that (\ref{M<}) holds 
a.s. on the event $\lef\{ \sum_{t \ge 0}\cR_t =\8 \rig\}$.
In particular, the inclusion opposite to (\ref{SsupL}) holds true.
\eds
\end{proposition}
%%%%%%%%%%%%%
We will prove \Prop{locb} via the following observation for 
general exponential martingales, which may be of independent interest.

Let $(M_t)_{t \in \N}$ be a square-integrable martingale 
on a filtered probability space 
$(\W, \cF, P \; ; \; (\cF_t)_{t \in \N})$. We denote its 
predictable quadratic variation by:
$$
\lan M \ran_t =
\sum_{1 \le u \le t}P[(\D M_u)^2 |\cF_{u-1}] 
$$
Here, and in what follows, we write $\D a_t=a_t-a_{t-1}$ ($t \geq
1$) for a sequence $(a_t )_{t \in \N}$ (random or non-random).
%%%%%%%%
\Proposition{gloc}
%%%%%%%
Let $(Y_t)_{t \in \N}$ be a mean-zero 
square-integrable martingale on a 
filtered probability space $(\W, \cF, P \; ; \; (\cF_t)_{t \in \N})$ 
such that $-1 \le \D Y_t$ a.s. for all $t \in \N^*$ 
and let 
%%%%%%%%
%\bdnl{gU_t}
%Y_t=\lef( {X_t \over X_{t-1}}-1 \ri){\bf 1}\{ X_{t-1} >0\},\; \; t \in \N^*,
%\edn
%so that
%%%%%%%% 
\bdnl{gU_t2}
X_t=\prod_{s=1}^t(1+\D Y_s).
\edn
\bds
\item[a)] 
Suppose that 
\bdnl{cond:X>0}
\sup_{t \ge 1}P[(\D Y_t)^2 |\cF_{t-1}] \le c_1^2\; \; a.s.
\edn
for some constant $c_1 \in (0,\8)$.
 Then, 
\bdnl{gloc1}
\{  X_\8 >0 \} \supset 
S \cap 
\lef\{ \;  \lan Y \ran_\8 <\8 \; \rig\}
\; \; a.s.
\edn
where $S=\{ \mbox{$X_t >0$ for all $t \ge 0$}\}$.
\item[b)]
Suppose that there exists a constant $c_2 \in (0,\8)$ such that 
for all $t \in \N^*$:
\bdnl{gU^3}
Y_t \in L^3 (P) \; \; \mbox{and}\; \; 
P[(\D Y_t)^3 |\cF_{t-1}] \le c_2P[(\D Y_t)^2 |\cF_{t-1}]
\; \; \mbox{a.s.}
\edn
Then, for any $c_3 \in (0,\; {1 \over 4})$, 
\bdnl{gM<}
X_t \le  \exp \lef( -c_3 \lan Y \ran_t \ri) 
\; \; \mbox{for all large enough $t$'s}
\edn
a.s. on the event 
$\lef\{ \;  \lan Y \ran_\8 =\8 \; \rig\}$.
In particular, the inclusion opposite to (\ref{gloc1}) holds true.
\eds
%%%%%
\end{proposition}
%%%%%%%%%
\noindent {\bf Remark:} 
As will be seen from the proof, 
the following assumption works as well for  \Prop{gloc}b): 
there exist $q \in (2,\8)$ and $c_2 \in (0,\8)$ 
such that for all $t \in \N^*$:
$$
Y_t \in L^q (P) \; \; \mbox{and}\; \; 
P[|\D Y_t|^q |\cF_{t-1}] \le c_2^{q-2}P[(\D Y_t)^2 |\cF_{t-1}]
\; \; \mbox{a.s.}
$$
Although this condition may look better than (\ref{gU^3}) for $q <3$, 
(\ref{gU^3}) works more effectively for our application. The point is that 
(\ref{gU^3}) is written in terms of $(\D Y_t)^3$, 
rather than $|\D Y_t|^3$.

\vvs
We postpone the proof of \Prop{gloc} (section \ref{p:gloc}) 
to finish the proof of \Prop{locb}.
%%%%%%%%%%%%

\vs
\noindent {\it Proof of \Prop{locb}}: 
%%%%%%%%%%%%%%%%%%%%%%%%%%%%%%%%%%%%%%%%%%%%%%%
We apply \Prop{gloc} to $X_t=|\ov{N}_t|$. 
Then, it is easy to see that (\ref{gU_t2}) holds with: 
$$
\D Y_t=
 {1 \over |a|}\sum_{x,y \in \zd}\rh_{t-1}(x)A_{t,x,y}-1 
$$
Moreover, it was shown in the proof of \cite[Lemma 3.2.2]{Yo08b} that 
there are constants $c_i \in (0,\8)$ ($i=1,2$) such that:
\bds 
\item[1)] \hspace{1cm}
$P\lef[ (\D Y_t)^p | \cF_{t-1} \ri]  \le   c_1\cR_{t-1},\; \; p=2,3$ 
\item[2)] \hspace{1cm}
$P\lef[ (\D Y_t)^2 | \cF_{t-1} \ri] 
 \ge   c_2\cR_{t-1}$.
\eds
((\ref{covA})) is used only for 2)). Therefore, 
\Prop{gloc} immediately leads to \Thm{loc}. 
\hfill $\Box$
%%%%%%%
\subsection{Proof of \Prop{gloc}}
\label{p:gloc}
%%%%%%%%
Let $(M_t)_{t \in \N}$ be a square-integrable martingale 
defined on a filtered probability space.
In this paper, we will repeatedly exploit 
the following  well-known facts (e.g., \cite[pages 252--253]{Dur05}):
\bdmn
 \{ \lan M \ran_\8 <\8\} & \sub & \{ \mbox{$M_t$ converges
as $t \ra \8$}\}\; \; \mbox{a.s.} \label{252} \\
\{ \lan M \ran_\8 =\8\} & \sub & 
\lef\{ \lim_{t \ra \8}{M_t \over \lan M \ran_t}=0\ri\}
\; \; \mbox{a.s.} \label{253}
\edmn

\vvs
To prove \Prop{gloc}, we will use the following lemma, 
which is a 
generalization of the Borel-Cantelli lemma, and is also 
used in the proof of \Lem{VW} below.
%%%%%%%%%
\Lemma{gBC}
%%%%%%%
Let $(Z_t)_{t \in \N}$ be an integrable, adapted
process defined on a filtered probability space 
$(\W, \cF, P \; ; \; (\cF_t)_{t \in \N})$ and let:
$$
A_0=0, \; \; 
A_t=\sum_{1 \le s \le t}P[\D Z_s|\cF_{s-1}],\; \; t \in \N^*.
$$
\bds
\item[a)] 
Suppose that there exists a constant $c_1 \in (0,\8)$ such that:
\bdnl{cond:gBC}
\D Z_t-P[\D Z_t|\cF_{t-1}] \ge -c_1\; \; 
\mbox{a.s. for all $t \in \N^*$.} 
\edn
Then, 
\bdnl{gBC}
\{\dps \lim_{t \ra \8}Z_t=\8 \} 
=
\lef\{\dps \lim_{t \ra \8}Z_t=\8, \; \; 
\suplim_{t \ra \8}{ A_t\over Z_t} \ge 1\rig\}
 \sub \{\dps \sup_{t \ge 1}A_t=\8 \} 
\; \;  \mbox{a.s.}
\edn
%%%%%%%%%%
%\barray{lll} 
%\{\dps \sum_{t \ge 1}X_t=\8 \} 
%& \sub &
%\lef\{\dps \inflim_{t \ra \8}
%{\sum_{1 \le s \le t}P[X_s|\cF_{s-1}] \over \sum_{s \le t}X_s} \ge 1\rig\}
%\; \;  \mbox{a.s.} \\
%& \sub &
%\{\dps \sum_{t \ge 1}P[X_t|\cF_{t-1}]=\8 \} 
%\earray 
%\edn
%%%%%%%%%
\item[b)] 
Suppose that $\{ Z_t \}_{t \in \N} \sub L^2 (P)$ and that 
there exists a constant $c_2 \in (0,\8)$ such that:
\bdnl{cond:gBC2}
{\rm var}(\D Z_t|\cF_{t-1}) \le c_2 P[\D Z_t|\cF_{t-1}]
\; \; \mbox{a.s. for all $t \in \N^*$,} 
\edn
where 
${\rm var}(\D Z_t|\cF_{t-1})=P[(\D Z_t)^2|\cF_{t-1}]-P[\D Z_t|\cF_{t-1}]^2$.
Then,
\bdnl{gBC2}
\{\dps \lim_{t \ra \8}A_t=\8 \} 
=
\lef\{\dps \lim_{t \ra \8}A_t=\8, \; \; 
\lim_{t \ra \8}{ Z_t \over A_t} = 1\rig\}
 \sub \{\dps \lim_{t \ra \8}Z_t=\8 \}.
\; \;  \mbox{a.s.}
\edn
%%%%%%%%%%%
\eds
%%%%%%%%%%
\end{lemma}
%%%%%%%%%%%%%
Proof: a) It is enough to show that
\bds \item[1)] \hspace{1cm}
$\{\dps \lim_{t \ra \8}Z_t=\8 \} 
 \sub 
\lef\{\dps 
\suplim_{t \ra \8}{ A_t\over Z_t} \ge 1\rig\}$.
\eds
Define $M_t=Z_t-A_t$, so that $(M_\cdot)$ is a martingale whose 
increments are bounded below by $-c_1$. Then, 
it is standard (e.g. the proof of \cite[page 236, (3.1)]{Dur05}) that
\bds \item[2)] \hspace{1cm}
$P(C \cup D_-)=1$,
\eds
where
$$
C=\{ \mbox{$M_t$ converges as $t \ra \8$}\} \; \; \mbox{and}\; \; 
D_-=\{ \inf_{t \in \N}M_t=-\8 \}.
$$
Now, by writing 
$$
{ A_t\over Z_t}=1-{ M_t\over Z_t},
$$
1) follows immediately from 2).\\
b) It is enough to show that
\bds \item[3)] \hspace{1cm}
$\{\dps \lim_{t \ra \8}A_t=\8 \} 
 \sub 
\lef\{\dps 
\lim_{t \ra \8}{ Z_t\over A_t} = 1\rig\}$.
\eds
Here, $M_\cdot$ is square-integrable. Since
$$
\lef| {Z_t \over A_t}-1 \ri|
= \lef|{ M_t \over A_t}\ri|,
$$
we have
$$
\{\dps \lim_{t \ra \8}A_t=\8, \; \lan M \ran_\8 <\8 \} 
\st{\scriptsize (\ref{252})}{\sub} 
\lef\{\dps 
\lim_{t \ra \8}{ Z_t\over A_t} = 1\rig\}.
$$
On the other hand, on the event $\{ \lan M \ran_\8 =\8\}$, 
we have 
$$
\lef|{ M_t \over A_t}\ri|
\st{\scriptsize (\ref{cond:gBC2})}{\le}
c_2{ |M_t| \over  \lan M \ran_t} 
\st{\scriptsize (\ref{253})}{\lra}0 
\; \; \; \mbox{as $t \ra \8$}
$$
These prove 3).
\hfill $\Box$

\vvs
\noindent {\bf Remark:} Similarly as \Lem{gBC}a), we can show the 
following variant of \Lem{gBC}b).
 Suppose that there exists a constant $c_3 \in (0,\8)$ such that:
$$
\D Z_t-P[\D Z_t|\cF_{t-1}] \le c_3\; \; 
\mbox{a.s. for all $t \in \N^*$.} 
$$
Then, 
$$
\{\dps \lim_{t \ra \8}A_t=\8 \} 
=
\lef\{\dps \lim_{t \ra \8}A_t=\8, \; \; 
\suplim_{t \ra \8}{ Z_t\over A_t} \ge 1\rig\}
 \sub \{\dps \sup_{t \ge 1}Z_t=\8 \} 
\; \;  \mbox{a.s.}
$$
%%%%%%%%%%
\Lemma{gloc(b)}
%%%%%%%%
Let $(Y_t)_{t \in \N^*}$ be as in \Prop{gloc}b). 
Then,
\bdnl{gloc(b)}
\{ \lan Y \ran_\8=\8 \}
\sub \lef\{ \inflim_{t \ra \8}
{\sum_{s \le t}f(\D Y_s) 
\over \lan Y \ran_t } \ge 1\rig\}
\; \; \; \mbox{a.s.}
\edn
where $f(u)={u^2 \over 2+u}$, $u \ge -1$.
%%%%%%%
\end{lemma}
%%%%%%%%%%%
Proof 
We first prepare elementary estimates. 
Let $U$ be a r.v. such that $-1 \le U$ a.s. 
Since 
$0 \le f(u) \vee f(u)^2 \le u^2$, we have 
\bds
\item[1)] \hspace{1cm}
$P\lef[f(U) \vee f(U)^2\ri] \le P[U^2]$.
\eds
Suppose further that $P[U^3] \le cP[U^2]$. Then,
\bds 
\item[2)] \hspace{1cm}
${\dps  P[U^2] \vee P[f(U)^2] \le (2+c)P\lef[f(U)\ri] }$.
\eds
This can be seen as follows. We have 
\bdnn
P[U^2]^2&=&P\lef[{U \over \sqrt{2+U}}U\sqrt{2+U}\rig]^2
 \le   P\lef[f(U) \ri] 
 P\lef[U^2 (2+U)\ri] \\
 & = & P\lef[f(U) \ri]
(2 P[U^2]+P[U^3]) 
 \le  (2+c)P[f(U)]P[U^2],
\ednn
which proves $P[U^2] \le (2+c)P\lef[f(U)\ri]$. 
On the other hand, 
$$
P[f(U)^2] \st{\mbox{\scriptsize 1)}}{\le} P[U^2] \le (2+c)P[f(U)].
$$
By 1)--2) above, applied to $U=\D Y_t$ and the 
measure $P(\; \cdot\; |\cF_{t-1})$, we see that
\bds 
\item[3)] \hspace{1cm}
${\dps
D\st{\rm def}{=}\{ \lan Y \ran_\8 =\8 \}
 =\{ \sum_{s \ge 1}P[f(\D Y_s)|\cF_{s-1}]=\8 \}\; \; }$a.s.
\eds
We see from 2) that $Z_t=\sum_{s \le t}f(\D Y_s)$ 
satisfies (\ref{cond:gBC2}). Therefore, 
$$
D \st{\scriptsize 3),(\ref{gBC2})}{\sub} \lef\{ \lim_{t \ra \8}
{\sum_{s \le t}f (\D Y_s) \over \sum_{1 \le s \le t}
P[f(\D Y_s)|\cF_{s-1}]}=1\rig\}
\; \; \; \mbox{a.s.}
$$
Thus, (\ref{gloc(b)}) follows from this and 1).
\hfill $\Box$
%%%%%%%%%%%%%%%%%%%%%%%%%%%%%

\vvs
%%%%%%%
\noindent {\it Proof of \Prop{gloc}}:a) 
We will prove that
\bds
\item[1)] \hspace{1cm}
$ S \cap \{ \lan Y \ran_\8 <\8\}
\sub \{ \exp (-Y_\8)X_\8 >0\}\; \; $ a.s.
\eds
We get (\ref{gloc1}) from this and (\ref{252}). 
To prove 1), note that
\bds
\item[2)] \hspace{1cm}
${\dps 
\exp (-Y_t)X_t=\prod^t_{u=1}(1+\D Y_u)\exp (-\D Y_u)
}$ \eds
and that
\bds
\item[3)] \hspace{1cm}
${\dps 
0 \le 1- ( 1+\D Y_u )\exp ( -\D Y_u )
\le {e \over 2}(\D Y_u)^2,
}$ \eds
since $\D Y_u \ge -1$. 
By (\ref{cond:X>0}), $Z_t=\sum_{s \le t}(\D Y_s)^2$ 
satisfies (\ref{cond:gBC}).
Thus, we have by (\ref{gBC}) that
\bds
\item[4)] \hspace{1cm}
${\dps 
 \{ \lan Y \ran_\8 <\8\}
\sub\{ \sum_{u \ge 1}(\D Y_u)^2 <\8 \}\; \;}$a.s. \eds
Thus, we get 1) from 2)--4). \\
b) We have $(1+u)e^{-u} \le e^{-f(u)/4}$ for 
$u \ge -1$, where  $f(u)={u^2 \over 2+u}$. Thus, 
\bds
\item[5)] \hspace{1cm}
$( 1+\D Y_u )\exp ( -\D Y_u ) 
\le \exp (-f(\D Y_u)/4)\; \; $ for all $u \ge 1$.
\eds
Let $0<c_3<c_4 <{1 \over 4}$.
Then, for $t$ large enough, 
a.s. on the event $\{ \lan Y \ran_\8 =\8\}$, 
\bdnn
\prod^t_{u=1}(1+\D Y_u)\exp (-\D Y_u) 
 &\st{\mbox{\scriptsize 5)}}{\le}&  
\exp \lef( - \sum^t_{u=1}f(\D Y_u)/4\ri) 
\st{\mbox{\scriptsize (\ref{gloc(b)})}}{\le}
\exp \lef( -c_4\lan Y \ran_t\ri) \\
&\st{\mbox{\scriptsize (\ref{253})}}{\le}&
\exp \lef( -Y_t-c_3\lan Y \ran_t\ri),
\ednn
which, via 2), proves (\ref{gM<}).
\hfill $\Box$
%%%%%%%
\subsection{Proof of \Thm{loc}a)}
%%%%%%%
If $P(|\ovN_\8|>0)>0$, then, 
$$
\lef\{ \mbox{survival} \rig\}=
\{ |\ovN_\8|>0\}\; \; \;\mbox{a.s.}
$$
This can be seen easily by translating the argument in 
\cite[page 701, proof of ``Proposition"]{Gri83}. We see from 
this and \Prop{locb} that $\sum_{t \ge 0}\cR_t <\8$ a.s. 
on the event of survival, 
while $\sum_{t \ge 0}\cR_t <\8$ is obvious outside the event of survival.
\hfill $\Box$
%%%%%%%%%%%%%%%%%%%%
%%%%%%%%%%%%%%%
\SSC{Proofs of \Thm{sloc}  and \Thm{sloc2}} \label{p:sloc12}
%%%%%%%%%
%%%%%%%%%%%%%%
\subsection{The argument by P. Carmona and Y. Hu} \label{p:sloc1}
%%%%%%%%%%%%%%%
For  $f,g:\zd \ra [0,\8)$,  
we define their convolution $f*g$ by: 
$$
(f*g)(x)=\sum_{y \in \zd}f(x-y)g(y), \; \; \; x \in \zd.
$$
For the notational convenience, we also write $a(y)$ for $a_y$.
We define: 
$$
b_t = \underbrace{b*...*b}_t
, \; \;  t \in \N^*\; \; 
\mbox{where 
$b(x)={1 \over |a|^2}\sum_{y \in \zd}a(y)a(y-x)$}, \\
%%%%%%%
%G_s &=&\del_0 +\sum_{u=1}^sb_s, \; \; t \in \N^* \cup \{ \8\}.
%%%%%%
$$
To interpret this, let $(\tl{S}_t)_{t \in \N}$ be the independent copy of 
$((S_t)_{t \in \N}, P^0_S)$, cf.(\ref{S_t}). 
Then, 
$$
b_t(x)=P^0_S \otimes P^0_{\tl{S}}(S_t -\tl{S}_t=x)
$$
Therefore, by (\ref{irred})
\bdnl{Sumr_t}
1+\sum_{t \ge 1}b_t (0) ={1 \over 1-\pi_d}
\lef\{ \barray{ll}
=\8 & \mbox{if $d =1,2$} \\
<\8 & \mbox{if $d \ge 3$}
\earray \rig.
\edn
We first note that there are $\e >0$ and $t_0 \in \N$ such that:
\bdnl{Sumr_t2}
\sum_{1 \le t \le t_0}b_t (0) \ge {1+\e \over \gm-1}.
\edn
For $d=1,2$, we take $\e =1$. Then, (\ref{Sumr_t2}) 
holds for $t_0$ large enough, since $\sum_{t \ge 1}b_t (0)= \8$. For
$d \ge 3$, the assumption (\ref{P(S=S)}) 
and (\ref{Sumr_t}) imply (\ref{Sumr_t2})
for small enough $\e >0$ and large enough $t_0$. 
We now fix $\e>0$  and $t_0$ and
define: 
\bdnl{X_t}
X_t = \lan g * \rh_t, \rh_t \ran,
\; \; \; \mbox{where $g = \sum_{s=1}^{t_0}b_s$.}
\edn
(The bracket $\lan \cdot, \; \cdot \ran$ stands for the 
inner product of $\ell^2 (\zd)$.)
Note that $0 \le g \in \ell^1 (\zd)$ and that
\bdnl{|X_t|}
|X_t| \le |(g * \rh_t)^2|^{1/2}|\rh_t^2|^{1/2} \le |g|\cR_t.
\edn
(Recall again that $|f|=\sum_x |f(x)|$ for $f:\zd \ra \R$).
Let: 
$$
X_t=M_t+A_t
$$
be Doob's decomposition, defined by:
\bdnl{A_t}
A_0=0, \; \;  \D A_t =P[\D X_t |\cF_{t-1}] 
\; \; \mbox{for $t \in \N^*$}.
\edn
Proof of \Thm{sloc} is based on the following two lemmas. 
%%%%%%%%%
\Lemma{s,t-s}
%%%%%%%%%%
There are constants $c_1,c_2 \in (0,\8)$ such that: 
$$
A_t \ge c_1 \sum_{0 \le u \le t-1} \cR_u
-c_2\sum_{0 \le u \le t-1} \cR_u^{3/2}
\; \; \; \mbox{for all $t \in \N^*$.}
$$
%%%%%%%%%%%%
\end{lemma}
%%%%%%%%%%%%%
%%%%%%%%%
\Lemma{VW}
%%%%%%%%%%
$$
\{ \sum_{u \ge 0} \cR_u =\8\} \sub
\lef\{ \lim_{t \ra \8}{M_t \over \sum_{0 \le u \le t-1} \cR_u}=0 \ri\}
\; \; \; \mbox{a.s.}
$$
%%%%%%%%%%%%
\end{lemma}
%%%%%%%%%%%%%
\noindent{\bf Proof of \Thm{sloc}:}
We may focus on the event 
$D=\{ \sum_{u \ge 0} \cR_u =\8\}$. 
It follows from  (\ref{|X_t|}) and \Lem{VW} that 
$$
\lim_{t \ra \8}{A_t \over \sum_{0 \le u \le t-1} \cR_u}=0
\; \; \; \mbox{a.s. on $D$}
$$
and hence from \Lem{s,t-s} that
$$
\inflim_{t \ra \8}{\sum_{0 \le u \le t-1} \cR_u^{3/2} 
\over \sum_{0 \le u \le t-1} \cR_u} \ge {c_1 \over c_2}
\; \; \; \mbox{a.s. on $D$}.
$$
This, together with (\ref{S=L}), proves \Thm{sloc}.
\hfill $\Box$
%%%%%%%%%
\subsection{Proof of \Lem{s,t-s}} \label{sec:s,t-s}
%%%%%%%%%%
The following technical lemma is an extension of 
\cite[Lemma 3.1.1]{HuYo08} to the case where the random variables
$U_i \ge 0$ may vanish with positive probability.
%%%%%%%%%%%
\Lemma{CH}
%%%%%%%%%
Let $U_i \ge 0$, $1 \le i \le n$ ($n \ge 2$) be independent
random variables such that: 
$$
\mbox{$P[U_i^{3} ] <\8$ for $i=1,..,n$ and 
$\sum_{i=1}^nm_i=1$,}
$$
where $m_i=P[U_i]$.Then, with $U=\sum_{i=1}^nU_i$,
\bdmn 
P\left[ {U_1U_2 \over U^2}: U >0 \right] 
& \ge & m_1m_2  -2m_2\mbox{var}(U_1)
-2m_1\mbox{var}(U_2),
\label{eta12} \\
P\left[{U_1^2 \over U^2}
: U >0 \right] & \ge &P[U_1^{2}]  \lef( 1+2m_1 \ri) -2P[U_1^{3}].
\label{eta11} \edmn 
%%%%%%
\end{lemma}
%%%%%%%%%%
Proof: 
Note that $x^{-2} \ge 3-2x$ for $x \in 
(0,\8)$. Thus, we have that 
\bdnn 
P\lef[ { U_1U_2 \over U^2 }:U>0 \ri] 
& \ge & P\lef[U_1U_2(3 -2U) :U>0 \ri]
=P\lef[U_1U_2(3 -2U) \ri] \\
& = & P\lef[U_1U_2(1 -2(U-1)) \ri] 
=m_1m_2-2 P\lef[U_1U_2(U-1) \ri],  \\
P\lef[U_1U_2(U-1) \ri] &= &P\lef[U_1U_2(U_1-m_1)\ri]
+P\lef[U_1U_2(U_2-m_2) \ri] \\
& = & m_2 \mbox{var}(U_1)+m_1\mbox{var}(U_2). 
\ednn 
These prove (\ref{eta12}). Similarly,
\bdnn 
P\lef[ { U_1^{2} \over U^2} : U>0\ri] 
& \ge & P\lef[U_1^{2}(3 -2U) : U>0\ri]
 =  P\lef[ U_1^{2}(3 -2U) \ri] \\
& = & P\lef[U_1^{2} \ri]-2P\lef[U_1^2(U-1) \ri], \\
P\lef[U_1^2(U-1) \ri] & = &P\lef[U_1^{2}(U_1-m_1) \ri]
=P\lef[U_1^{3} \ri]-m_1P\lef[U_1^{2} \ri]. 
\ednn These prove
(\ref{eta11}). \hfill $\Box$

\vvs
We introduce 
\bdnl{rht,s} 
\rh_{t,1}=\rh_t * \ov{a}, 
\; \; \; \cR_{t,1}=|\rh_{t,1}^2|,
\edn
where $\ov{a}(x)=a(x)/|a|, \; \; x \in \zd$.

We will make a series of estimates on quantities involving
$a (x)$, $\rh_t (x)$, $\cR_t$, and so on. In the sequel, 
multiplicative constants are denoted by $c, c_1, c_2,..$.
We agree that they are {\it non-random} constants 
which do not 
depend on time variables  $t,s,.. \in \N$ 
or space variables $x,y,...\in \zd$.
%%%%%%%%%%%%%%
\Lemma{<cR<}
%%%%%%%%%
For any $t \in \N$,
\bdnl{<cR<}
\cR_{t,1} \le  \cR_t \le {|a|^2 \over |a^2|}\cR_{t,1}.
\edn
%%%%%%%%%%%
\end{lemma}
%%%%%%%%
Proof: Let $\ov{a}(x)=a(x)/|a|, \; \; x \in \zd$. We then have 
$$
|\rh_{t,1}^2|=|(\rh_t*\ov{a})^2| \le 
|\rh_t^2|
$$
by Young's inequality. This proves the first inequality. 
On the other hand, 
\bdnn
|\rh_{t,1}^2| &= & |(\rh_t*\ov{a})^2|
=\sum_{x \in \zd}
\lef( \sum_{y \in \zd}\rh_{t}(x-y)\ov{a}(y)\ri)^2 \\
& \ge & 
\sum_{x \in \zd}\sum_{y \in \zd}\rh_{t}(x-y)^2\ov{a}(y)^2
=|\rh_t^2||\ov{a}^2|,
\ednn 
which proves the second inequality. 
\hfill $\Box$

\vvs
We assume (\ref{covA}) from here on.
%%%%%%%%%%%%%%%%
\Lemma{P[rhrh]}
%%%%%%%%%%%
There is a constant $c \in (0,\8)$ such that the following hold:
\bdmn 
\lefteqn{P\lef[ \rh_{t} (y)\rh_{t} (\tl{y}) | \cF_{t-1} \ri]} \nn \\
& \ge &
\rh_{t-1,1} (y)\rh_{t-1,1} (\tl{y})
-c\rh_{t-1,1} (y)\rh_{t-1,1} (\tl{y})^2
-c\rh_{t-1,1} (\tl{y})\rh_{t-1,1} (y)^2,\label{ynot=tly}
\edmn
for all $t \in \N^*$, $y, \tl{y} \in \zd$ with $y \neq \tl{y}$.
\bdnl{y=tly}
P\lef[ \cR_t  | \cF_{t-1} \ri] 
\ge \gm \cR_{t-1,1}-c\cR_{t-1,1}^{3/2}
\; \; \; \mbox{for all $t \in \N^*$}.
\edn
%%%%%%%%%
\end{lemma}
%%%%%%%%%
Proof: Let $U_t=\sum_{y \in \zd}U_{t,y}$, where 
$U_{t,y}={1 \over |a|}\sum_{x \in \zd}\rh_{t-1}(x)A_{t,x,y}$. 
Then, $\{U_{t,y}\}_{y \in \zd}$ 
are independent under $P(\cdot |\cF_{t-1})$.
Moreover, it is not difficult to see that 
(cf. proof of \cite[Lemma 3.2.2]{Yo08b}), on the event 
$\{ |\ovN_{t-1}|>0\}$, 
\bds
\item[1)] \hspace{1cm}
$P[U_{t,y}|\cF_{t-1}] = \rh_{t-1,1}(y), \; \; \; P[U_t|\cF_{t-1}] = 1$,
\item[2)] \hspace{1cm}
${\dps P[U_{t,y}^2|\cF_{t-1}] ={1 \over |a|^2}
\sum_{x_1,x_2,y \in \zd}\rh_{t-1}(x_1)\rh_{t-1}(x_2)
P[A_{t,x_1,y}A_{t,x_2,y}]}$ 
\item[3)] \hspace{1cm}
${P[U_{t,y}^m|\cF_{t-1}]  \le  c_1  \rh_{t,1}(y)^m,\; \; \; m=2,3.}$
\eds
%%%%%%%
%where the constant $c_3 \in (0,\8)$ is independent of 
%$t$ and $y$ . 
%%%%%%%%%
Since 
$$
\rh_{t} (y)\rh_{t} (\tl{y})=(U_{t,y}U_{t,\tl{y}}/U_t)
{\bf 1}_{\{|\ovN_{t-1}|>0\}}
$$
and $\{ U_t>0 \} \sub \{|\ovN_{t-1}|>0\}$, we see from 
1), 3) above and \Lem{CH} that (\ref{ynot=tly}) holds and that
\bds
\item[4)] \hspace{1cm}
${\dps P\lef[ \rh_{t} (y)^2 | \cF_{t-1} \ri] 
 \ge
P[U_{t,y}^2|\cF_{t-1}]-2c_1\rh_{t-1,1} (y)^3}$.
\eds
To prove (\ref{y=tly}), note that 
\bds
\item[5)] \hspace{1cm}
${\dps \sum_{y \in \zd}\rh_{t-1,1} (y)^3 
\le \lef( \sum_{y \in \zd}\rh_{t-1,1} (y)^2\ri)^{3/2} =\cR_{t-1,1}^{3/2}}$.
\eds
We then see that
\bdnn
P\lef[ \cR_t | \cF_{t-1} \ri] 
& \st{\scriptstyle 4)}{\ge} & 
\sum_{y \in \zd}\lef( P[U_{t,y}^2|\cF_{t-1}]-2c_1\rh_{t-1,1} (y)^3\ri) \\
& \st{\scriptstyle 2), 5)}{\ge} &  
{1 \over |a|^2} \sum_{x_1,x_2,y \in \zd}
\rh_{t-1}(x_1)\rh_{t-1}(x_2)P[A_{t,x_1,y}A_{t,x_2,y}] 
-2c_1\cR_{t-1,1}^{3/2}\\
& \st{\scriptstyle (\ref{covA})}{\ge} & 
 { \gm\over |a|^2} \sum_{x_1,x_2,y \in \zd}
\rh_{t-1}(x_1)\rh_{t-1}(x_2)a(y-x_1)a(y-x_2)-2c_1\cR_{t-1,1}^{3/2} \\
& = & \gm \cR_{t-1,1} -2c_1\cR_{t-1,1}^{3/2}.
\ednn
\hfill $\Box$
%%%%%%%%%%%%%%%%%%%

\vvs
\noindent{\bf Proof of \Lem{s,t-s}}:
%%%%%%
$$
P[X_t|\cF_{t-1}]=
\sum_{y,\tl{y}\in \zd}g(y-\tl{y})
P[\rh_t (y)\rh_t (\tl{y})|\cF_{t-1}]=I+J,
$$
where $I$ and $J$ are diagonal and off-diagonal terms:
\bdnn
I &= &g(0)\sum_{y\in \zd}P[\rh_t (y)^2|\cF_{t-1}], \\
J &= &\sum_{y,\tl{y}\in \zd \atop y \neq \tl{y}}
g(y-\tl{y})P[\rh_t (y)\rh_t (\tl{y})|\cF_{t-1}].
\ednn
We start with the lower bound for $I$.
%%%%%%%%%
\bds \item[1)] \hspace{1cm}${\dps 
I=   g(0)P[\cR_t|\cF_{t-1}] 
 \st{\scriptstyle (\ref{y=tly})}{\ge} 
g(0)\gm \cR_{t-1,1}-g(0)c\cR_{t-1,1}^{3/2}. 
}$ \eds
%%%%%%%%%%%%%%
As for $J$, we have 
$$
J\st{\scriptstyle (\ref{ynot=tly})}{\ge} 
J_{1,1}  - cJ_{1,2}- cJ_{2,1},
$$
where
$$
J_{m,n}
= \sum_{y,\tl{y}\in \zd \atop y \neq \tl{y}}
g(y-\tl{y})\rh_{t-1,1} (y)^m\rh_{t-1,1} (\tl{y})^n.
$$
$J_{1,1}$ can be computed exactly:
%%%%%%%%%
\bds \item[2)] \hspace{1cm}${\dps 
\barray{ccl}
J_{1,1}
& = & \lef(\sum_{y,\tl{y}\in \zd}
-\sum_{y,\tl{y}\in \zd \atop y = \tl{y}}\ri)
g(y-\tl{y}) \rh_{t-1,1} (y)\rh_{t-1,1} (\tl{y}) \\
&= & \lan g * b* \rh_{t-1}, \rh_{t-1} \ran -g(0)\cR_{t-1,1}.
\earray 
}$ \eds
%%%%%%%%%%%%%%
To bound $J_{1,2}$ from above, note that
$$
\max_{x \in \zd}(g*\rh_{t-1,1})(x) 
\le |g|\max_{x \in \zd}\rh_{t-1,1}(x) \le |g|\cR_{t-1,1}^{1/2}.
$$
Thus, 
$$
J_{1,2}
 \le 
\lan g * \rh_{t-1,1}, \rh_{t-1,1}^2 \ran
\le \max_{x \in \zd}(g*\rh_{t-1,1})(x) \cR_{1,t-1}
\le |g|\cR_{t-1,1}^{3/2}.
$$
Similarly, $J_{2,1} \le |g|\cR_{t-1,1}^{3/2}$. Putting things together, 
we see that
\bds \item[3)] \hspace{1cm}${\dps 
\barray{ccl}
\D A_t & = & P[X_t|\cF_{t-1}]-X_{t-1} 
 \ge  I+J_{1,1}-X_{t-1} -2c|g|\cR_{t-1,1}^{3/2} \\
& \st{\mbox{\scriptsize 1)--2)}}{\ge} &
 (\gm -1)g(0)\cR_{t-1,1}
+\lan (g * b-g)* \rh_{t-1}, \rh_{t-1} \ran -3|g|c\cR_{t-1,1}^{3/2}.
\earray }$
\eds
Note that $g * b-g =b_{t_0+1}-b \ge -b$ and hence that
\bds \item[4)] \hspace{1cm}${\dps 
\lan (g * b-g)* \rh_{t-1}, \rh_{t-1} \ran
\ge -\lan b* \rh_{t-1}, \rh_{t-1} \ran
=-\cR_{t-1,1}.}$
\eds
Therefore,
$$
\D A_t
 \st{\mbox{\scriptsize 3)--4)}}{\ge} 
((\gm -1)g(0)-1)\cR_{t-1,1} -3|g|c\cR_{t-1,1}^{3/2} 
 \st{\mbox{\scriptsize (\ref{Sumr_t2})}}{\ge} 
 \e \cR_{t-1,1}-3|g|c\cR_{t-1,1}^{3/2}.
$$
We now get \Lem{s,t-s} from this and (\ref{<cR<}).
\hfill $\Box$
%%%%%%%%%
\subsection{Proof of \Lem{VW}} \label{sec:VW}
%%%%%%%%%%
We have 
$$
\lef\{ \sum_{1 \le u < \8}\cR_u =\8, \; \lan M \ran_\8 <\8 \ri\}  
\st{\scriptsize (\ref{252})}{\sub}
\lef\{ 
\lim_{t \ra \8}{M_t 
\over \sum_{1 \le u \le t}\cR_u }=0 \rig\}
\; \; \; \mbox{a.s.}
$$
To treat the case of $\lan M \ran_\8 =\8$, we show that
\bds \item[1)] \hspace{1cm}
${\dps \lan M \ran_t 
\le 4|g|^2\sum_{1 \le u \le t}\lef( \cR_{u-1}+P[ \cR_u |\cF_{u-1}]\ri) }$.
\eds
We have 
\bds \item[2)] \hspace{1cm}
$|\D X_t|^2  \st{\mbox{\scriptsize (\ref{|X_t|})}}{\le} 
2|g|^2 (\cR_t^2+\cR_{t-1}^2) 
\le 2|g|^2 (\cR_t+\cR_{t-1})$,
\eds
and 
\bds \item[3)] \hspace{1cm}
$(\D A_t)^2 
\st{\mbox{\scriptsize Schwarz}}{\le} 
P[(\D X_t)^2 |\cF_{t-1}] 
\st{\mbox{\scriptsize 2)}}{\le} 
 2|g|^2 (P[\cR_t |\cF_{t-1}] +\cR_{t-1}).$
\eds
Thus, 
\bdnn 
\D \lan M \ran_t & = & 
P[ (\D M_t)^2 |\cF_{t-1}]  
\le  2P[(\D X_t)^2|\cF_{t-1}] +2(\D A_t)^2  \\
&\st{\mbox{\scriptsize 3)}}{\le} &
4|g|^2 (P[\cR_t |\cF_{t-1}] +\cR_{t-1}).
\ednn
Now, we have by \Lem{gBC} and 1) that
\bdnn
\lef\{ \sum_{1 \le u < \8}\cR_u =\8 \ri\}  
&\st{\mbox{\scriptsize (\ref{gBC})}}{=} &
\lef\{ \sum_{1 \le u < \8}P[ \cR_u |\cF_{u-1}]=\8 \rig\}
\; \; \mbox{a.s.} \\
&\st{\mbox{\scriptsize (\ref{gBC2})}}{=} &
\lef\{ 
\lim_{t \ra \8}{\sum_{1 \le u \le t}\cR_u 
\over \sum_{1 \le u \le t}P[ \cR_u |\cF_{u-1}]}=1 \rig\} 
\; \; \mbox{a.s.}\\
&\st{\mbox{\scriptsize 1)}}{\sub} &
\lef\{ 
\inflim_{t \ra \8}{\sum_{1 \le u \le t}\cR_u
\over  \lan M \ran _t } \ge {1 \over 4|g|^2} \rig\}
\ednn
We see from this and (\ref{253}) that
$$
\lef\{ \sum_{1 \le u < \8}\cR_u =\8, \; \lan M \ran_\8 =\8 \ri\}  
\sub
\lef\{ 
\lim_{t \ra \8}{M_t 
\over \sum_{1 \le u \le t}\cR_u }=0 \rig\}
\; \; \; \mbox{a.s.}
$$
This completes the proof of \Lem{VW}.
\hfill $\Box$
%%%%%%%%%
\subsection{Proof of \Thm{sloc2}} \label{p:sloc2}
%%%%%%%%%%%%%%%%%%%%%
We now state a criterion for the regular growth phase (\Lem{Birk}). 
The criterion is an extension of the one obtained 
by M. Birkner \cite{Bir04} for DPRE. 

Let $((S_t)_{t \in \N}, P_S^x)$ be the random walk defined by 
(\ref{S_t}) 
and let $(\tl{S}_t)_{t \in \N}$ be its independent copy. 
Since the random variable:
$$
V_\8(S, \tl{S})=\sum_{t \ge 1}{\bf 1}_{\{ S_t=\tl{S}_t\}}
$$
is geometrically distributed with the parameter $\pi_d$, we have 
\bdnl{1/pi}
{1 \over \pi_d}=\sup\lef\{ \a  \ge 1 \; ; \; 
P_S^0 \otimes P_{\tl{S}}^0\lef[ \a^{V_\8(S, \tl{S})}\ri] <\8 \ri\}.
\edn
We now define $\pi^*_d$ by:
\bdnl{1/pi^*}
{1 \over \pi^*_d}=\sup\lef\{ \a  \ge 1 \; ; \; 
P_{\tl{S}}^0\lef[ \a^{V_\8(S, \tl{S})}\ri] <\8 \; \; 
\mbox{$P^0_S$-a.s.} \ri\}.
\edn
Therefore, $\pi^*_d \le \pi_d $ in general. Moreover, the inequality 
is known to be strict if $d \ge 3$ and (\ref{Cor4}) is satisfied 
\cite[page 82, Corollary 4]{Bir03}.
%%%%%%%%%%
\Lemma{Birk}
%%%%%%
Suppose $d \ge 3$ and (\ref{ovh1}).
Then, 
$$
P[\ov{\h}_{t,y}^2] < {1 \over \pi^*_d} \; \Ra \; P [|\ovn_\8|]=1.
$$
%%%%%%%%%%
\end{lemma}
%%%%%%%%%
Proof: Because of (\ref{ovh1}), we have that 
$$
N_{t,x}=|a|^tP_S^0 \lef[ \prod^t_{u=1}\ov{\h}_{u,S_u}\ri].
$$
Using this expression, we can repeat 
the argument in \cite{Bir04} without change.
(Here, unlike the DPRE case, we may have $P(\ov{\h}_{t,y}=0)>0$. 
However, this does not cause any problem as far as to prove 
this lemma.) \hfill $\Box$
%%%%%%%%%%%%%%%%%%5

\vvs 
\noindent {\bf Proof of \Thm{sloc2}:}
%%%%%%%%%%%
(\ref{limR>c}) $\sub$: 
Note that $\pi^*_d < \pi_d $ if $d \ge 3$ and (\ref{Cor4}) is satisfied.
If $|\ovn_\8|=0$ a.s., then we have by 
\Lem{Birk} that $\gm \ge {1 \over \pi^*_d}> {1 \over \pi_d} $. 
Thus, we can apply \Thm{loc} and \Thm{sloc}. \\
%%%%%%%%%%%%%%%%%%%%%%%%%%%%%%%
(\ref{limR>c}) $\supset$: Obvious.
\hfill $\Box$

\vvs
\footnotesize
%%%%%%%%%%%%%
\noindent{\bf Acknowledgements:}
The author thanks Yukio Nagahata for discussions which lead to 
the simplification of the proof of \Thm{sloc}. 

%%%%%%%%%%%%%%

%%%%%
%%%%%%%%%%
\end{document}